# A Better Upper Bound of Hausdorff Measure of the Cartesian Product of the Middle Third Cantor Set with Itself Compare to Others'

## FAN Yu-chen

(The Affiliated High School of South China Normal University, 510630, Guangzhou, Guangdong, China)

**Abstract:** After the correction of an inaccurate result in the reference, the author uses five different methods, and gets five different inequalities on the Hausdorff measure of the Cartesian product of the middle third Cantor set with itself:

$$\begin{cases} H^s(C \times C) \leq 1.548563 \\ H^s(C \times C) \leq 1.504975 \\ H^s(C \times C) \leq 1.502878 \\ H^s(C \times C) \leq 1.503263 \\ H^s(C \times C) \leq 1.502483 \end{cases}$$

The main theory of this paper is the best among them, which is

$$H^s(C \times C) \leq 1.502483.$$

**Key Words:** Self-similar Set; The Cartesian Product of the Middle Third Cantor Set with Itself, Hausdorff Measure and Dimension; Partial Estimation Principle.

## 1 Introduction

The calculation of Hausdorff measure for fractal sets is one of the most important subjects in fractal geometry, and it is a very difficult research topic [1]. Since then, only a small number of Hausdorff measures of fractal sets with Hausdorff dimension smaller than 1 have been computed, for example, the Cantor and Sierpinski blanket [2,3], and, no Hausdorff measure of a fractal set with Hausdorff dimension greater than 1 has been calculated [4]. In this paper, the author gives three different methods and use these methods to obtain better estimste on the upper bound of the Hausdorff measure of the Cartesian product of the middle third Cantor set with itself.

Cutting out a unit square out of the plane, retaining 4 squares with side length 1/3 on the 4 corners of the original square with side length 1, then take out the rest of the original square, and repeating these operations to the 4 squares on the corners of the original square for infinity number of times, we obtain a set of points which is called "the Cartesian product of the middle third Cantor set with itself", and is denoted by "$C \times C$" [5]. Both Cantor and $C \times C$ are uncountable sets (see 4.1).

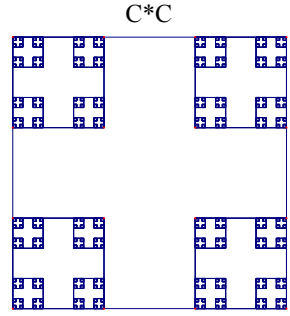

C*C

In this paper, the 1st level subset is the "biggest" or the "first" subset, for example, the 1st level subset, or say the 1st level small square is the square with side length 1.

## 2 Main Result

**Theorem: The Hausdorff measure of the Cartesian product of the middle third Cantor set with itself satisfies the inequality**
$$H^s(C \times C) \leq 1.502483.$$

## 3 Definitions

**3.1 Definition 1** [1]**: Self-similar Set:** Assume $E \subset R^n$ to be a bounded-closed subset, and $f: E \to E$ to be a mapping. If $\forall x, y \in E, \exists c > 0$ s.t.
$$|f(x) - f(y)| \leq c|x - y|,$$
then we call $f$ is a "Compression Mapping" or "Compression Function" on $E$ with compression ratio $c$. Obviously, compression functions are continuing. When the equality holds, we call $f$ is a "Similar Compression Function" on $E$ with similarity ratio $c$. Let $\{f_1, f_2, \ldots, f_n\}$ be a Compression Function System on $E$ with similarity ratio $0 < c_i < 1, \ i = 1,2, \ldots, m$. There exists a unique subset $E_s$ of $E$ that satisfies
$$E_s = \bigcup_{i=1}^{n} f_i(E_s),$$
then we call $E_s$ be an "Invariant Set" of $\{f_1, f_2, \ldots, f_n\}$ [1]. For any $p > 1$, assume the invariant set of $\{f_1, f_2, \ldots, f_n\}$ to be
$$E_s = \bigcup_{i=1}^{n} f_i(E_s) = \bigcup_{J_p} f_{i_1} f_{i_2} \cdots f_{i_k}(E_s).$$
For $(i_0 i_1 \ldots i_{k-1}) \in J_p$, let $f_i(E_s)$ be the first iteration of $E_s$, where $i = 0,1, \ldots, \infty$; inductively, let $f_{i_0} f_{i_1} \ldots f_{i_{k-1}}(E_s)$ be the $k^{\text{th}}$ iteration of set $E_s$. Obviously, every iteration of $E_s$ is similar to $E_s$. This is the definition of self-similar set. In short, after enlarging arbitrary part of a self-similar set with proper ratio, this part will coincide exactly with the original set.

The construction of the middle third Cantor set: define a similar compression function system $\{f_1, f_2\}$ on the interval $E = [0,1]$:
$$f_1(x) = \frac{x}{3}, \qquad f_2(x) = \frac{x+2}{3}.$$

The invariant set $E_s$ we get on the interval $E$ is the middle third Cantor set.

**3.2 Definition 2** [5]**: Hausdorff measure:** Let $E$ be a subset of $R^n$, if there exist nonempty subsets $U_i$, $i = 1,2,...$ of set $R^n$ such that

$$E \subset \bigcup_{i=1}^{\infty} U_i, 0 < |U_i| \leq \delta.$$

Then this is called a "$\delta$ cover" of set $E$. For $s \geq 0$, $\delta > 0$, we use "inf" to represent taking the infimum of $\delta$ cover of the set $E$, and define

$$H_\delta^s(E) = \inf \sum_{i=1}^{\infty} |U_i|^s \geq 0.$$

Pick $0 < \delta_1 < \delta_2$, then we have $0 < |U_i| \leq \delta_1$ and $0 < |U_i| \leq \delta_2$. Because the infimum of the integral set must be smaller than the infimum of the partial set, it is going to be

$$H_{\delta_1}^s(E) \geq H_{\delta_2}^s(E).$$

As a result, function $f(\delta) = H_\delta^s(E)$ is increasing in its domain. When $\delta \to 0^+$,

$$H^s(E) = \lim_{\delta \to 0^+} H_\delta^s(E) \geq 0$$

and $H^s(E)$ is called the Hausdorff measure of $E$ in $s$-dimension. Hausdorff measure is used to distinguish the "size" between countable sets and uncountable sets.

**3.3 Definition 3** [5]**: Hausdorff dimension:** For arbitrary set $E$, when $s$ is increasing in the interval $(0,\infty)$, $H_\delta^s(E)$ does not increase. When $s < j$, the inequality

$$H_\delta^s(E) \geq \delta^{s-j} H_\delta^j(E)$$

always hold. That is, if $H_\delta^j(E)$ is positive, $H_\delta^s(E)$ tends to infinite. As a result, there is exactly a critical value $s$ that satisfies

$$\begin{cases} H^s(E) = \infty \ (0 \leq s < \dim_H(E)) \\ H^s(E) = 0 \ (\dim_H(E) < s < \infty) \end{cases},$$

such $s = \dim_H(E)$ is called the Hausdorff dimension of $E$.

Hausdorff dimension of $C \times C$: $C \times C$ is made up of 4 similar compression functions with similar rate 1/3, and it satisfies the condition of opened set. So, the Hausdorff dimension of $C \times C$ satisfies

$$4 \cdot \left(\frac{1}{3}\right)^s = 1.$$

That is,

$$s = \dim_H(C \times C) = \log_3 4.$$

## 4 Lemmas

**4.1 Lemma 1:** The middle third Cantor set is an uncountable set.

**Proof:** Transfer all real numbers in interval $(0,1)$ into base-3 form, then the middle third Cantor set can be represented in this form:

$$C = \begin{cases} 0.a_{11}a_{12}a_{13}\ldots a_{1n}\ldots \,|\, a_{1i} \in \{0,2\}, i = 1,2,\ldots,n,\ldots \\ 0.a_{21}a_{22}a_{23}\ldots a_{2n}\ldots \,|\, a_{2i} \in \{0,2\}, i = 1,2,\ldots,n,\ldots \\ \ldots \\ 0.a_{k1}a_{k2}a_{k3}\ldots a_{kn}\ldots \,|\, a_i \in \{0,2\}, i = 1,2,\ldots,n,\ldots \end{cases}$$

for all position integer $k$. Pick

$$\varepsilon = 0.a_1^* a_2^* \ldots a_n^* \ldots$$

satisfying $a_j^* \neq a_{jj}, j = 1,2,\ldots,n,\ldots, a_j^* \in \{0,2\}$. According to the construction of $\varepsilon$, we know that $\varepsilon \in C$, but $\varepsilon \notin C$ according to the construction, contradiction! As a result, it is impossible to list all elements in Cantor, Cantor is "unlistable", that is, it is uncountable. According to Lemma 1, it is obviously that $C \times C$ is an uncountable set.

**4.2.1 Lemma 2:** Partial Estimation Principle: Construct measurable set $U \subset R^n$, the diameter of $U$ is $|U|$, and $|U| > 0$, then the inequality

$$H^s(E \cap U) \leq |U|^s.$$

(in this inequality, $E$ is the set to be measured and $E \subset R^n$ is a Self-similar set satisfying the condition of open set; $s$ is the Hausdorff dimension of $E$) holds [1, 8].

**Proof: Sub-lemma 1:** Construct $S_1 = \{U_l | l \geq 0\}$ such that it is $\delta$ cover of $E$. Let $S_2 = \{U_i | i \geq 0\}$ be another $\delta$ cover of $E$ ($S_1 = S_2$ is not necessary needed). By the definition, we have

$$H_\delta^s(E) = \inf \sum_{i=1}^{\infty} |U_i|^s \leq \sum_{l=0}^{\infty} |U_l|^s.$$

Assume $P_{j_1}, \ldots, P_{j_k}(\alpha)$ is a $\delta_k = (\max_j\{c_j\})^k \delta$ cover of $E$, then we have

$$H_{\delta_k}^s(S_1) \leq \sum_{J_k}\sum_{l=0}^{\infty} c_{j_1}^s \ldots c_{j_k}^s |U_l|^s = \sum_{l=0}^{\infty}\left(\sum_{j_1}^{m} c_{j_1}^s\right)\ldots\left(\sum_{j_k}^{m} c_{j_k}^s\right)|U_l|^s = \sum_{l=0}^{\infty} |U_l|^s.$$

When $k \to \infty$, $\delta_k \to 0$, we have

$$H^s(E) \leq \sum_{l=0}^{\infty} |U_l|^s.$$

**Sub-lemma 2:** Because the Hausdorff dimension of $E - E \cap U$ is not larger than $s$, we have $H^s(E - E \cap U) < \infty$. By the definition of Hausdorff measure, for any $\varepsilon > 0$, there exists a cover $S_3 = \{U_i' | i > 0\}$ of set $E - E \cap U$ satisfying

$$H^s(E - E \cap U) + \varepsilon \geq \sum_{i=0}^{\infty} |U_i'|^s.$$

**Proof of Lemma 2:** According to sub-lemma 1 and sub-lemma 2, because $U$ is a cover set of $E$ and $E \cap U$, cover $S_4 = \{U_k | k > 0\} = \{U_i', U | i > 0\}$ is a cover of $E$, we must have

$$H^s(E) = H^s(E \cap U) + H^s(E - E \cap U) \leq \sum_{i=0}^{\infty} |U_k|^s = \sum_{i=0}^{\infty} |U_i'|^s + |U|^s,$$

$$\sum_{i=0}^{\infty} |U_i'|^s + |U|^s \leq H^s(E - E \cap U) + |U|^s + \varepsilon.$$

So, we have

$$H^s(E \cap U) \leq |U|^s + \varepsilon,$$

let $\varepsilon \to 0$, then we have the result

$$H^s(E \cap U) \leq |U|^s.$$

**4.2.2 Inference 1:** For arbitrary positive integer $k$, assume the number of $k^{\text{th}}$ level elements of $E$ (a self-similar open set with similar ratio $c$) to be $m^{k-1}$, and the number of the $k^{\text{th}}$ level elements cover by cover set to be $N$. Then we have

$$\frac{N}{m^{k-1}} H^s(E) \leq |U|^s.$$

**Proof:** $E$ is made up of $m$ similar compression functions with similar ratio $c$, so the Hausdorff dimension of this set $s = \dim_H(E)$ satisfies $c^s = 1/m$. Pick cover $U$ as

$$\bigcup_{i=1}^{N} S_i^k \subset U;$$

where $S_i^k$, $i = 1, 2, \ldots, N$ is the $k^{\text{th}}$ level subset of $E$. According to the partial estimation principle, we have

$$|U|^s \geq H^s(E \cap U) \geq H^s\left(\bigcup_{i=1}^{N} S_i^k\right) = Nc^{s \cdot (k-1)} H^s(E) = \frac{N}{m^{k-1}} H^s(E).$$

As a result, we have

$$\frac{N}{m^{k-1}} H^s(E) \leq |U|^s.$$

**4.2.2 Inference 2:** For our measured set $E$ (a self-similar open set), we have

$$H^s(E) \leq |E|^s$$

**Proof:** Pick cover set $U = E$, according to Inference 1, we have

$$\frac{E}{E} H^s(E) = H^s(E) \leq |E|^s.$$

As a result, the inequality

$$H^s(E) \leq |E|^s$$

holds.

**4.3 Lemma 3:** Among all sets with the same diameter in the same Euclidean plane, the circle has the biggest area [9].

**Proof:** Assume the diameter of a set $U$ is $d$ with $d > 0$. Take an arbitrary point $O$ on the edge of set $U$. Construct a polar coordinate system with point $O$ as the original point. Assume the polar equation of bound of $U$ be $R(\theta)$, where $0 \leq \theta \leq \pi$. Pick points

$$(R(\theta), \theta), \quad \left(R\left(\theta + \frac{\pi}{2}\right), \theta + \frac{\pi}{2}\right).$$

According to the Pythagorean theorem

$$d_0 = \sqrt{R(\theta)^2 + R\left(\theta + \frac{\pi}{2}\right)^2},$$

notice that $d_0 \leq d$, that is,

$$d_0^2 = R(\theta)^2 + R\left(\theta + \frac{\pi}{2}\right)^2 \leq d^2.$$

Therefore, the area bounded by the edge of $U$ is

$$S_U = \int_b^a \frac{R(\theta)^2}{2} d\theta = \int_0^\pi \frac{R(\theta)^2}{2} d\theta =$$

$$\int_{\frac{\pi}{2}}^\pi \frac{R(\theta)^2}{2} d\theta + \int_0^{\frac{\pi}{2}} \frac{R(\theta)^2}{2} d\theta = \int_0^{\frac{\pi}{2}} \frac{R\left(\theta + \frac{\pi}{2}\right)^2 + R(\theta)^2}{2} d\theta \leq$$

$$\int_0^{\frac{\pi}{2}} \frac{d^2}{2} d\theta = \frac{d^2}{2} \cdot \frac{\pi}{2} - \frac{d^2}{2} \cdot 0 = \frac{\pi d^2}{4}.$$

Notice that the area of circle with diameter $d$ is

$$S_0 = \pi r^2 = \pi \left(\frac{d}{2}\right)^2 = \frac{\pi d^2}{4} \geq S_U.$$

As a result, the circle has the biggest area among all sets with the same diameter in the same Euclidean plane. Lemma 4 explains why we use circles as the cover set in the proof of the theorem below.

## 5 Correction

The correction of the inaccurate upper bound in reference [10]: The author of reference [10] calculated an upper bound of the Hausdorff measure of the Cartesian product of the middle third Cantor set with itself:

$$H^s(C \times C) \leq 1.495901.$$

However, the author of this paper found that this upper bound is inaccurate and corrected this upper bound. After the correction, the new upper bound is:

$$\boldsymbol{H^s(C \times C) \leq 1.512163}.$$

**Proof:** According to the method of construct $U$ in reference [10], construct the rectangular coordinator plane $xA_1y$ with point $A_1(0,0)$, then construct square $A_1A_2A_3A_4$ with $A_1A_2=1$ by point $A_1$ for vertex and $x$, $y$ axes for sides, take points

$$B_1(2/3,0), B_2(2/3,9), B_3(25/3,9), B_4(25/3,0).$$

Use $B_1$ as the center, segment $B_1B_3$ as the radius, and construct an arc intersecting line $x=9$ in the first quadrant at $C_1$. Use $B_2$ as the center, segment $B_2B_4$ as the radius, construct an arc intersecting line $x=9$ in the first quadrant at $C_2$. Use $B_3$ as the center, segment $B_1B_3$ as the radius, and construct an arc intersects the $y$ axis in the first quadrant at $C_3$. Use $B_4$ as the center, segment $B_2B_4$ as the radius, and construct an arc intersecting the $y$ axis in the first quadrant at $C_4$. Take the area bounded by arcs $B_1C_3$, $B_2C_4$, $B_3C_1$ and $B_4C_2$ and segments $B_1B_4$, $B_2B_3$, $C_1C_2$ and $C_3C_4$ as the cover set $U$. Notice that the construction of $U$ is more complicated because there are 4 different equations of circles. However, consider the symmetry and the rotation of this figure, we only need to calculate the arc constructed by center $B_3$ and segment $B_1B_3$. The whole figure of this method is shown below.

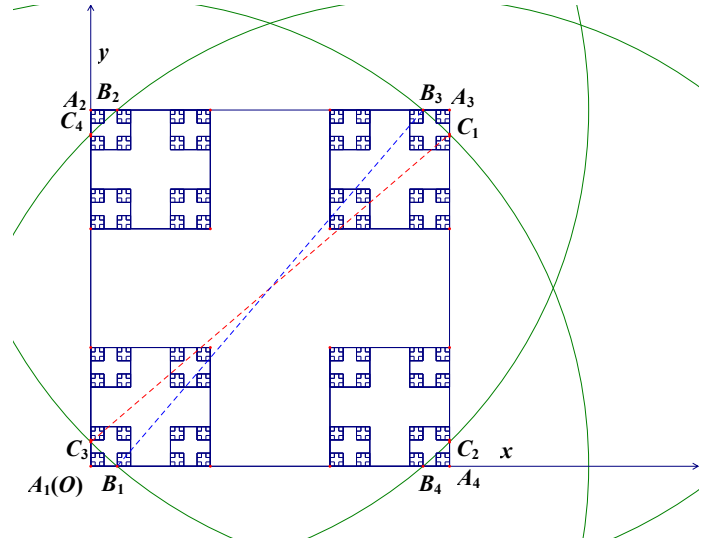

The radius of this circle is
$$r = \sqrt{\left(9-\frac{4}{3}\right)^2 + 9^2} = \frac{\sqrt{1258}}{3}.$$

So, the equation of this circle is
$$\left(x-\frac{25}{3}\right)^2 + (y-9)^2 = \frac{1258}{9}.$$

Taking
$$x = \frac{1}{3},$$
substituting $x$ into the equation of $\odot B_3$ and picking the smaller $y$, we have
$$y = 9 - \frac{\sqrt{682}}{3} \approx 0.2950 < \frac{1}{3} \approx 0.3333.$$

Taking
$$x = 0,$$
substituting $x$ into the equation of $\odot B_3$ and picking the smaller $y$, we have
$$y = 9 - \frac{\sqrt{633}}{3} \approx 0.6135 < \frac{2}{3} \approx 0.6667.$$

According to the symmetry of this figure, we know that $\odot B_3$ must intersect with the 4[th] level small square lying on the bottom left corner, but does not intersect with the 4[th] level small square lying on the top and near of this small square. **The 1[st] level square here is the square $A_1A_2A_3A_4$ with side length 9.**

Taking
$$x = \frac{2}{9},$$
substituting $x$ into the equation of $\odot B_3$ and picking the smaller $y$, we have
$$y = 9 - \frac{\sqrt{5993}}{9} \approx 0.3984 > \frac{2}{9} \approx 0.2222.$$

Taking

$$y = \frac{1}{3},$$

substituting $x$ into the equation of $\odot B_3$ and picking the smaller $x$, we have

$$x = \frac{25}{3} - \frac{\sqrt{582}}{3} \approx 0.2918 > \frac{2}{9} \approx 0.2222.$$

The calculation above shows that, when

$$x = \frac{1}{3},$$

the value of $y$ is

$$y = 9 - \frac{\sqrt{682}}{3} \approx 0.2950 < \frac{1}{3} \approx 0.3333.$$

Now, we know that $\odot B_3$ must pass through the 5$^{th}$ level small square on the top right corner of the 4$^{th}$ level small square lying on the bottom left corner, as shown in the figure below.

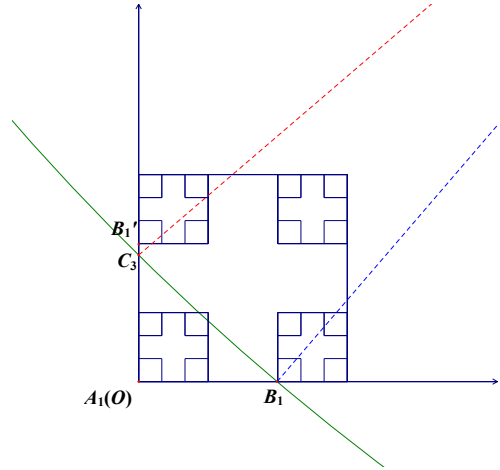

Consider the coverage of 6$^{th}$ level squares.
The calculation above shows that, when

$$y = \frac{1}{3},$$

the value of $x$ is

$$x = \frac{25}{3} - \frac{\sqrt{582}}{3} \approx 0.2918 > \frac{2}{9} + \frac{1}{27} \approx 0.2593,$$

$$x = \frac{25}{3} - \frac{\sqrt{582}}{3} \approx 0.2918 < \frac{2}{9} + \frac{2}{27} \approx 0.2963.$$

When

$$x = \frac{1}{3},$$

the value of $y$ is

$$y = 9 - \frac{\sqrt{682}}{3} \approx 0.2950 > \frac{2}{9} + \frac{1}{27} \approx 0.2593,$$

$$y = 9 - \frac{\sqrt{682}}{3} \approx 0.2950 < \frac{2}{9} + \frac{2}{27} \approx 0.2963.$$

Taking

$$x = \frac{2}{9} + \frac{2}{27} = \frac{8}{27},$$

substituting $x$ into the equation of $\odot B_3$ and picking the smaller $y$, we have

$$y = 9 - \frac{\sqrt{54809}}{27} \approx 0.3291 > \frac{2}{9} + \frac{2}{27} \approx 0.2963.$$

Now, we know that $\odot B_3$ must pass through the $6^{th}$ level small square on the top right corner of the $5^{th}$ level small square on the top right corner of the $4^{th}$ level small square lying on the bottom left corner, as shown in the figure below.

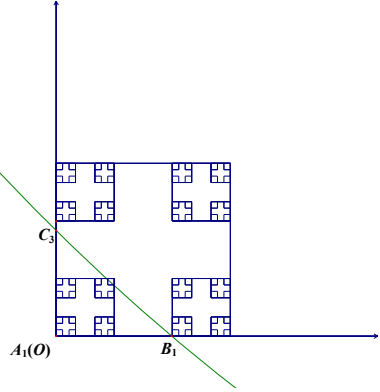

Consider the coverage of $7^{th}$ level squares.
The calculation above shows that, when

$$x = \frac{2}{9} + \frac{2}{27} = \frac{8}{27},$$

the value of $y$ is

$$y = 9 - \frac{\sqrt{54809}}{27} \approx 0.3291 > \frac{2}{9} + \frac{2}{27} + \frac{2}{81} \approx 0.3210,$$

$$y = 9 - \frac{\sqrt{54809}}{27} \approx 0.3291 < \frac{1}{3} \approx 0.3333.$$

Taking

$$x = \frac{2}{9} + \frac{2}{27} + \frac{2}{81} = \frac{26}{81},$$

substituting $x$ into the equation of $\odot B_3$ and picking the smaller $y$, we have

$$y = 9 - \frac{\sqrt{495881}}{81} \approx 0.3063 > \frac{2}{9} + \frac{2}{27} \approx 0.2963,$$

$$y = 9 - \frac{\sqrt{495881}}{81} \approx 0.3063 < \frac{2}{9} + \frac{2}{27} + \frac{1}{81} \approx 0.3086.$$

Taking

$$x = \frac{2}{9} + \frac{2}{27} + \frac{1}{81} = \frac{25}{81},$$

substituting $x$ into the equation of $\odot B_3$ and picking the smaller $y$, we have

$$y = 9 - \frac{\sqrt{494582}}{81} \approx 0.3177 < \frac{1}{3} - \frac{1}{81} \approx 0.3210,$$

$$y = 9 - \frac{\sqrt{495881}}{81} \approx 0.3063 < \frac{2}{9} + \frac{2}{27} + \frac{1}{81} \approx 0.3086.$$

The calculation above shows that, when

$$x = \frac{1}{3},$$

the value of $y$ is

$$y = 9 - \frac{\sqrt{682}}{3} \approx 0.2950 < \frac{2}{9} + \frac{2}{27} \approx 0.2963.$$

Now, we know that $U$ must pass through three $7^{th}$ level small squares on the top right corner of the $6^{th}$ level small square on the top right corner of the $5^{th}$ level small square on the top right corner of the $4^{th}$ level small square lying on the bottom left corner, as shown in the figure below. Inside, $U$ fully contains the $7^{th}$ level small squares on the top right corner of this $6^{th}$ level small square, and $U$ also contains another two $7^{th}$ level small squares, as shown in the figure below.

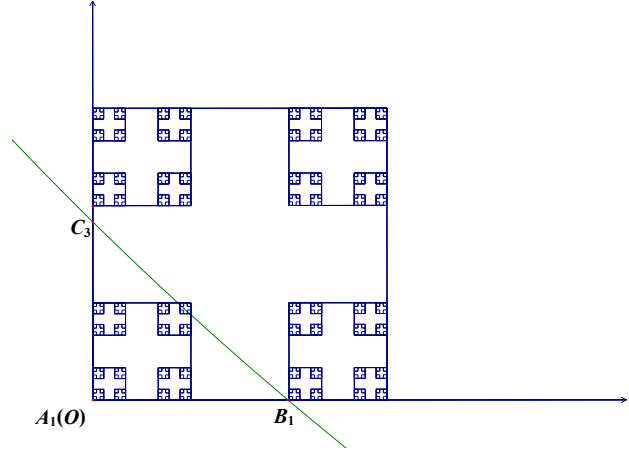

Consider the coverage of $8^{th}$ level squares.
Taking

$$y = \frac{1}{3} - \frac{1}{27} = \frac{8}{27},$$

substituting $x$ into the equation of $\odot B_3$ and picking the smaller $y$, we have

$$x = \frac{25}{3} - \frac{\sqrt{46673}}{27} \approx 0.3319 > \frac{1}{3} - \frac{1}{243} \approx 0.3292.$$

Taking

$$x = \frac{2}{9} + \frac{2}{27} + \frac{1}{243} = \frac{73}{243},$$

substituting $x$ into the equation of $\odot B_3$ and picking the smaller $y$, we have

$$y = 9 - \frac{\sqrt{4443434}}{243} \approx 0.3253 > \frac{1}{3} - \frac{1}{243} \approx 0.3251.$$

Taking

$$y = \frac{1}{3} - \frac{1}{81} = \frac{26}{81},$$

substituting $x$ into the equation of $\odot B_3$ and picking the smaller $x$, we have

$$x = \frac{25}{3} - \frac{\sqrt{422873}}{81} \approx 0.3051 > \frac{1}{3} - \frac{2}{81} - \frac{1}{243} \approx 0.3045.$$

The calculation above shows that, when

$$x = \frac{2}{9} + \frac{2}{27} = \frac{8}{27},$$

the value of $y$ is

$$y = 9 - \frac{\sqrt{54809}}{27} \approx 0.3291 < \frac{1}{3} - \frac{1}{243} \approx 0.3292,$$

$$y = 9 - \frac{\sqrt{54809}}{27} \approx 0.3291 > \frac{1}{3} - \frac{2}{243} \approx 0.3251.$$

The calculation above shows that, when

$$x = \frac{2}{9} + \frac{2}{27} + \frac{2}{81} = \frac{26}{81},$$

the value of $y$ is

$$y = 9 - \frac{\sqrt{495881}}{81} \approx 0.3063 > \frac{1}{3} - \frac{2}{81} - \frac{1}{243} \approx 0.3045.$$

Taking

$$x = \frac{1}{3} - \frac{2}{243} = \frac{79}{243},$$

substituting $x$ into the equation of $\odot B_3$ and picking the smaller $y$, we have

$$y = 9 - \frac{\sqrt{4466822}}{243} \approx 0.3025 < \frac{1}{3} - \frac{2}{81} - \frac{2}{243} \approx 0.3004.$$

This result claims that, in these three $7^{th}$ level small squares above, there are almost three $8^{th}$ level small squares on two $7^{th}$ level small squares on the top left corner and bottom right corner cover by $U$, respectively. Also, four $8^{th}$ level small squares on the $7^{th}$ level small squares on the top right corner are covered by $U$ too, as shown in the figure below.

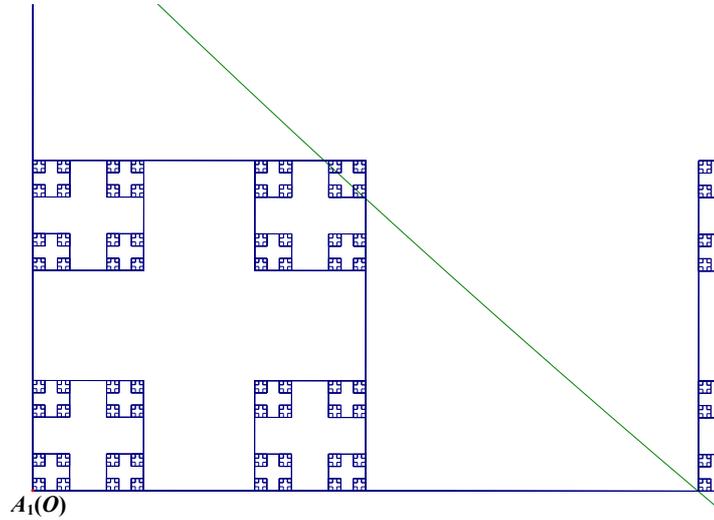

According to the symmetry and rotation of $U$, the number of small squares that do not cover by $U$ is

$$4 \times (4^2 \times 15 + 6) = 984,$$

and the total number of $8^{th}$ level small squares is

$$4^7 = 16384.$$

As a result, the fraction of $8^{th}$ level small squares contained by $U$ is

$$\frac{4^7 - 984}{4^7} = \frac{1925}{2048}.$$

Now we calculate the diameter of $U$.
Taking
$$x = 0,$$
substituting $x$ into the equation of $\odot B_3$ and picking the smaller $y$, we have
$$y = 9 - \frac{\sqrt{633}}{3} \approx 0.6135 < \frac{2}{3} \approx 0.6667.$$
This result shows that actually the coordinator of point $C_1$ is under (0,2/3), according to the symmetry, *actually, the diameter is $C_1C_3$.*
Proof: $C_4B_2$ and $B_4C_2$ below are two arcs in the figure shown below. Take point $K \in C_4B_2$, if
$$\exists M \in B_4C_2 \text{ s.t.} \max\{|KM|\},$$
we must have $M = C_2$ or $M = B_4$. When $M = C_2$, if
$$\exists K \in C_4B_2 \text{ s.t.} \max\{|KM|\},$$
we must have $K = C_4$; when $M = B_4$, if
$$\exists K \in C_4B_2 \text{ s.t.} \max\{|KM|\},$$
we must have $K = B_2$. According to the article above, we know that
$$C_4C_2 > B_4B_2.$$
As a result, the diameter of $U$ is $C_4C_2$, according to the symmetry, the diameter is $C_1C_3$, which is the red dotted segment below.

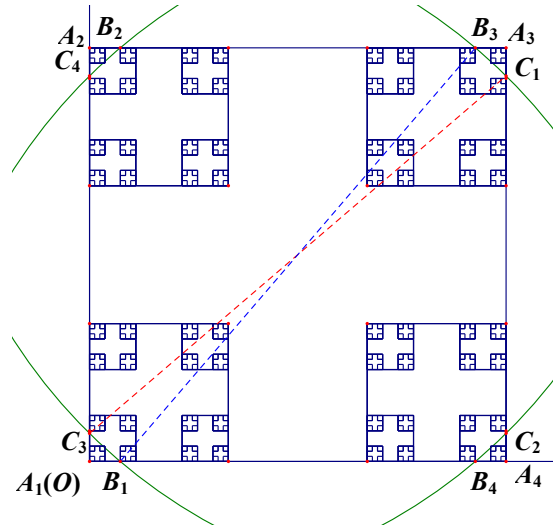

As a result,
$$|U| = \frac{\sqrt{\left(\frac{2\sqrt{633}}{3} - 9\right)^2 + 9^2}}{9} \approx 1.32133.$$
So, we have
$$\frac{1925}{2048} H^s(C \times C) \leq |U|^s = (1.32133)^{\log_3 4}.$$
Solving this inequity, we have
$$\boldsymbol{H^s(C \times C) \leq 1.512163}.$$
Under the situation of this method, the upper bound of
$$H^s(C \times C)$$

is bigger than this result, because some 8th level small squares are not fully contained by $U$ (With this method, the author assume they are fully contained).

The mistake made by the author of reference [10] is that the author didn't find the diameter correctly. However, the author of this article successfully finds the correct diameter of set $U$ and gets an accurate result.

## 6 Theorems

### 6.1 Attentions

**1.** When the side length of the 1st level square is $m$, $\forall k \in N^+, k \geq 1$, and the side length of $k^{th}$ level small square is
$$m \cdot 3^{1-k}.$$
**2.** According to the definition, the side length of the biggest square of the Cartesian product of middle third Cantor set with itself should be 1. However, in order to calculate it more easily, the side length of the biggest square of Cartesian product of middle third Cantor set with itself is changed into 9. In order to get the right answer, we only need to divide the diameter by 9, and this action will not influence the result of this paper.

### 6.2 Theorem 1

Can we just pick a subset of $C \times C$ as the cover set? Notice that the area of one $n^{th}$ level small squares of $C \times C$ is $(1/3)^{2(n-1)}$, and the number of $n^{th}$ level squares is $4^{n-1}$, so the cover fraction is $4^{-(n-1)}$, where the 1st level square is the square with length 1 here. Take one square from $n^{th}$ level squares as the cover set $U$, then the diameter of $U$ is equal to $\sqrt{2} \cdot (1/3)^{(n-1)}$. By Lemma 2, the Hausdorff dimension of $C \times C$ is $s = \dim_H(C \times C) = \log_3 4$. As a result, according to the partial estimation principle, we have

$$4^{-(n-1)} H^s(C \times C) \leq \left[\sqrt{2} \cdot \left(\frac{1}{3}\right)^{n-1}\right]^{\log_3 4}.$$

Solving this inequality, we have

$$H^s(C \times C) \leq \sqrt{2}^{\log_3 4} \approx 1.548563.$$

When $n = 2$, the cover set $U$ with its diameter $|U|$ are shown in the figure below.

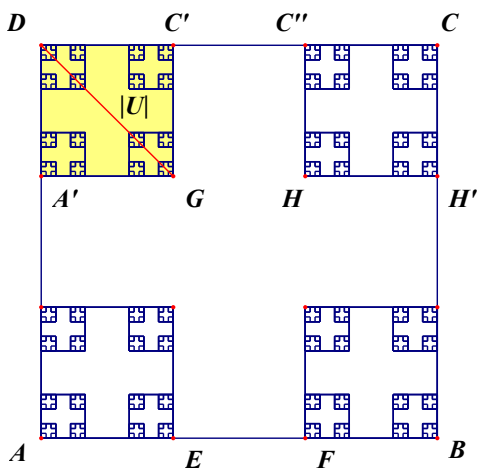

This kind of cover set is called the "Basic Interval" of a fractal set. However, the upper bound above shows that we should find another cover set to measure or estimate the Hausdorff measure of $C \times C$ because it is obviously not good enough.

**6.2 Theorem 2**

In order to compute the upper bound of
$$H^s(C \times C)$$
correctly, we can take the cover set $U$ as an octagon. In this method, the Hausdorff measure of the Cartesian product of the middle third Cantor set with itself satisfies the inequality
$$\boldsymbol{H^s(C \times C) \leq 1.504975}.$$

**Proof:** Firstly, construct the rectangular coordinate system $xA_1y$ with point $A_1(0,0)$, then construct square $A_1A_2A_3A_4$ with $A_1A_2=1$ by point $A_1$ for vertex and $x$, $y$ axis for sides. Notice that for arbitrary non-negative integer $n$, $n+1$ is the "level number" of the little square(s) (that is, 1st level, 2nd level, …, n+1th level). For $C \times C$, the number of elements it contains is $4^n$, the area of every element is $1/9^n$, and the side length of every element is $1/3^n$. Take points

$B_1$ (2/27, 9), $B_2$ (2/27, 9), $B_3$ (25/27, 9), $B_4$ (25/27, 0);
$C_1$ (0, 2/27), $C_2$ (0, 25/27), $C_3$ (9, 25/27), $C_4$ (9, 2/27),

and construct an octagon $B_1C_1C_2B_2B_3C_3C_4B_4$ as cover set $U$ in lemma 3, as shown in the figure below.

**The 1st level square here is the square $A_1A_2A_3A_4$ with side length 1.**

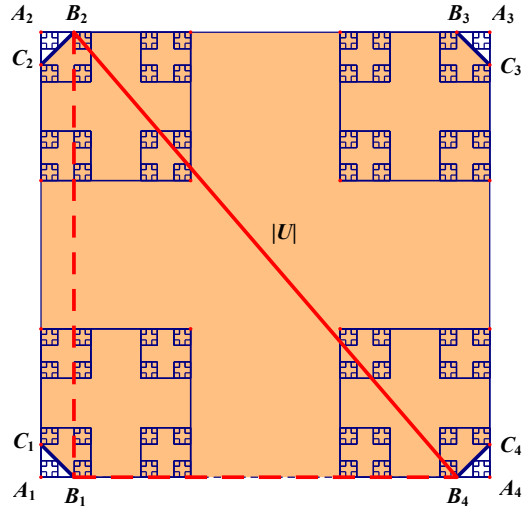

Notice that the diameter of this cover set is
$$|U| = \sqrt{B_1B_2^2 + B_1B_2^2} = \sqrt{1 + \left(\frac{23}{27}\right)^2} = \frac{\sqrt{1258}}{27}.$$

The fraction of squares cover by $U$ is
$$1 - \frac{4}{4^3} = \frac{15}{16}.$$

Then by Lemma 3, we have
$$\frac{15}{16} H^s(C \times C) \leq \left(\frac{\sqrt{1258}}{27}\right)^{\log_3 4}.$$

Solving this inequality, we have
$$H^s(C \times C) \leq 1.504975.$$
However, this upper bound is not optimized enough. We will use another method.

**6.3 Theorem 2**

We also use an octagon as cover set $U$ in this method. However, the difference is, that this method uses constructing infinity series and estimating the value of
$$H^s(C \times C).$$
With this method, the Hausdorff measure of the Cartesian product of the middle third Cantor set with itself satisfies the inequality
$$H^s(C \times C) \leq 1.502878.$$

**Proof:** Firstly, construct the rectangular coordinate system $xA_1y$ with point $A_1(0,0)$, then construct square $A_1A_2A_3A_4$ with $A_1A_2=1$ by point $A_1$ for vertex and $x$, $y$ axes for sides. Notice that for arbitrary non-negative integer $n$, $n+1$ is the "level number" of the little square(s) (that is, 1st level, 2nd level, …, n+1th level). For $C \times C$, the number of elements it contains is $4^n$, the area of every element is $1/9^n$, and the side length of every element is $1/3^n$. With this method, the author constructs an octagon inside this square as cover $U$ in lemma 3, as shown in the figure below.

**The 1st level square here is the square $A_1A_2A_3A_4$ with side length 1.**

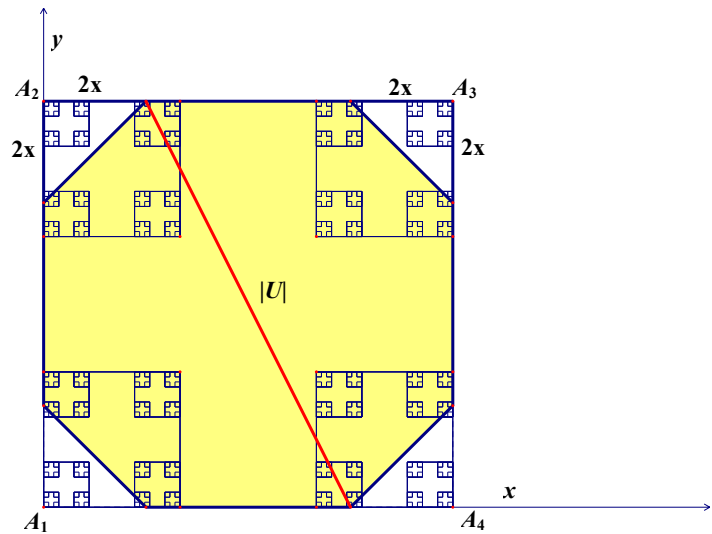

Construct this octagon such that for arbitrary $n$, the octagon doesn't include all-level small squares in the corner of the big square. Taking
$$x = \sum_{i=1}^{\infty} \frac{1}{3^{ki}}$$
for arbitrary positive integer $k$, the geometrical meaning of $x$ here is the sum of side lengths of all $(jk+1)$th level small squares, where $j$ is a positive integer. By the symmetry of the octagon, the diameter of $U$ is
$$|U| = \sqrt{1 + (1-4x)^2}.$$
As it is known from the geometrical meaning and symmetry, when $j$ plus 1, the total number of small square that *do not covered by U* on the top left corner is increased by A Better Upper Bound of Hausdorff Measure of the Cartesian Product of the Middle Third Cantor Set with Itself Compare to Others'$2^{i-1}$, and the total number of small squares in this level is increased by $4^{jk}$. As a result, the ratio of squares covered by $U$ is

$$N = 1 - 4\sum_{i=1}^{\infty} \frac{2^{i-1}}{4^{ki}}.$$

According to the figure, we must have

$$1 - 4x \geq 0, x \leq \frac{1}{4}.$$

When $k = 1$, the value of $x$ is

$$x = \sum_{i=1}^{\infty} \frac{1}{3^i} = \frac{1}{2} > \frac{1}{4}.$$

As a result, $k \neq 1$. Because $k$ is a positive integer, we have $k \geq 2$. Assume

$$k = 2 \Rightarrow x_1 = \sum_{i=1}^{\infty} \frac{1}{3^{2i}} = \frac{1}{8},$$

$$k = 3 \Rightarrow x_2 = \sum_{i=1}^{\infty} \frac{1}{3^{3i}} = \frac{1}{26},$$

$$k = 4 \Rightarrow x_3 = \sum_{i=1}^{\infty} \frac{1}{3^{4i}} = \frac{1}{80},$$

$$k = 5 \Rightarrow x_4 = \sum_{i=1}^{\infty} \frac{1}{3^{5i}} = \frac{1}{242};$$

then the frictions of squares covered by $U$ for each different values of $x$ are

$$x_1 = \frac{1}{8} \Rightarrow N_1 = 1 - 4\sum_{i=1}^{\infty} \frac{2^{i-1}}{4^{2i}} = \frac{5}{7},$$

$$x_2 = \frac{1}{26} \Rightarrow N_2 = 1 - 4\sum_{i=1}^{\infty} \frac{2^{i-1}}{4^{3i}} = \frac{29}{31},$$

$$x_3 = \frac{1}{80} \Rightarrow N_3 = 1 - 4\sum_{i=1}^{\infty} \frac{2^{i-1}}{4^{4i}} = \frac{125}{127},$$

$$x_4 = \frac{1}{242} \Rightarrow N_4 = 1 - 4\sum_{i=1}^{\infty} \frac{2^{i-1}}{4^{5i}} = \frac{509}{511};$$

then by Lemma 2, we have

$$\frac{5}{7} H^s(C \times C) \leq \left(\sqrt{1 + (1-4x)^2}\right)^s = \left(\frac{\sqrt{5}}{2}\right)^{\log_3 4},$$

$$\frac{29}{31} H^s(C \times C) \leq \left(\sqrt{1 + (1-4x)^2}\right)^s = \left(\frac{\sqrt{290}}{13}\right)^{\log_3 4},$$

$$\frac{125}{127} H^s(C \times C) \leq \left(\sqrt{1+(1-4x)^2}\right)^s = \left(\frac{\sqrt{761}}{20}\right)^{\log_3 4},$$

$$\frac{509}{511} H^s(C \times C) \leq \left(\sqrt{1+(1-4x)^2}\right)^s = \left(\frac{\sqrt{28802}}{121}\right)^{\log_3 4}.$$

Therefore, for each different values of $k$, we have
$$k = 2 \implies H^s(C \times C) \leq 1.611653,$$
$$k = 3 \implies H^s(C \times C) \leq 1.502878,$$
$$k = 4 \implies H^s(C \times C) \leq 1.524502,$$
$$k = 5 \implies H^s(C \times C) \leq 1.538520.$$

So, when
$$x = \sum_{i=1}^{\infty} \frac{1}{3^{3i}} = \frac{1}{26},$$

the upper bound of $H^s(C \times C)$ is the best upper bound, given by
$$\mathbf{H^s(C \times C) \leq 1.502878.}$$

However, mathematics is based on proof. In our observation we find that when $k = 3$, the upper bound is the best. However, we still need to proof this result. We can use $k$ to represent the upper bound of $H^s(C \times C)$, and find the value of $k$ satisfy the minimum value of the upper boun of $H^s(C \times C)$. For every positive integer $k$, we have

$$x = \sum_{i=1}^{\infty} \frac{1}{3^{ki}} = \frac{\frac{1}{3^k}}{1-\frac{1}{3^k}} = \frac{\frac{1}{3^k}}{\frac{3^k-1}{3^k}} = \frac{1}{3^k-1}.$$

then the frictions of squares covered by $U$, $N$, is

$$N = 1 - 4\sum_{i=1}^{\infty} \frac{2^{i-1}}{4^{ki}} = 1 - 4 \cdot \frac{\frac{1}{4^k}}{1-\frac{2}{4^k}} = 1 - 4 \cdot \frac{\frac{1}{4^k}}{\frac{4^k-2}{4^k}} = 1 - \frac{4}{4^k-2} = \frac{4^k-6}{4^k-2}.$$

Submit the values of $x$ and $N$ into Inference 1 of Lemma 2, we have

$$N \cdot H^s(C \times C) = \frac{4^k-6}{4^k-2} \cdot H^s(C \times C) \leq \left(\sqrt{1+(1-4x)^2}\right)^s$$

$$= \left(\sqrt{1+\left(1-\frac{4}{3^k-1}\right)^2}\right)^s = \left(\sqrt{1+\left(\frac{3^k-5}{3^k-1}\right)^2}\right)^s = \left(\sqrt{\frac{(3^k-1)^2+(3^k-5)^2}{(3^k-1)^2}}\right)^s$$

$$= \left(\frac{\sqrt{(3^k-1)^2+(3^k-5)^2}}{3^k-1}\right)^s = \left(\frac{\sqrt{3^{2k}-2 \cdot 3^k+1+3^{2k}-10 \cdot 3^k+25}}{3^k-1}\right)^s$$

$$= \left(\frac{\sqrt{2 \cdot 3^{2k}-12 \cdot 3^k+26}}{3^k-1}\right)^s.$$

As a result, we have

$$H^s(C \times C) \leq \left(\frac{\sqrt{2 \cdot 3^{2k} - 12 \cdot 3^k + 26}}{3^k - 1}\right)^s \cdot \frac{4^k - 2}{4^k - 6}.$$

Assume $k$ can be nonnegative real numbers, define function $f$ as

$$f(k) = \left(\frac{\sqrt{2 \cdot 3^{2k} - 12 \cdot 3^k + 26}}{3^k - 1}\right)^s \cdot \frac{4^k - 2}{4^k - 6},$$

Where $k$ defines at $[2, +\infty)$. Then we are going to take the first derivative for this function. Use Matlab (Matlab code see appendix) to find the first derivative $f'(k)$ of this function and solve $k$ for $f'(k) = 0$, we have

$$k_0 = 2.78051450628200301042367443964451.$$

Take the second derivative $f''(k)$ and submit this value of $k$ into it and find

$$f''(k_0) = 0.1063 > 0.$$

As a result, when $k = x_0$, the function achieves its minimum value. Consider back when $k$ is an integer greater than 1, we can pick $k = 2$ and $k = 3$, then we find that when $k = 3$, the upper bound of $H^s(C \times C)$ achieve its minimum value, that is, it is the best result in this method. The graph of $f(k)$ is shown below (Here we do not consider the graph when $k \in [0,2)$, as I explained before).

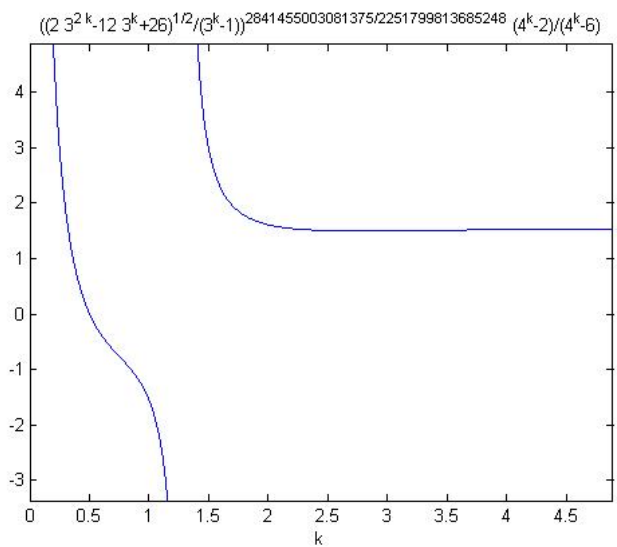

However, we can question that the upper bound obtained by this method may not be optimized enough. Hence the author will find another method, the idea is to expands the cover ratio of $U$ but keep the diameter of $U$. That is, the octagon may not be chosen as the cover set $U$. The cover set $U$ can be chosen from sets having same diameter as the octagons but having more cover ratio than the octagons (in short, the area is bigger). By Lemma 4, circle or arc can be chosen as cover set $U$.

**6.4 Theorem 3**

Consider to change the cover ratio of $U$, keep the diameter of $U$ stay constant, so that the upper bound of

$$H^s(C \times C).$$

can be smaller; the Hausdorff measure of the Cartesian product of the middle third Cantor set with itself satisfies the inequality

$$H^s(C \times C) \leq 1.503263.$$

(In this method, the author expands the diameter in sub theorem 1 but still uses circle as the cover set, because in this way the upper bound can be smaller).

**Proof:** Firstly, construct the rectangular coordinator system $xA_1y$, Point $A_1(0,0)$, then construct square $A_1A_2A_3A_4$ with $A_1A_2=9$ by point $A_1$ for vertex and $x$, $y$ axes for sides, as shown in the figure below.

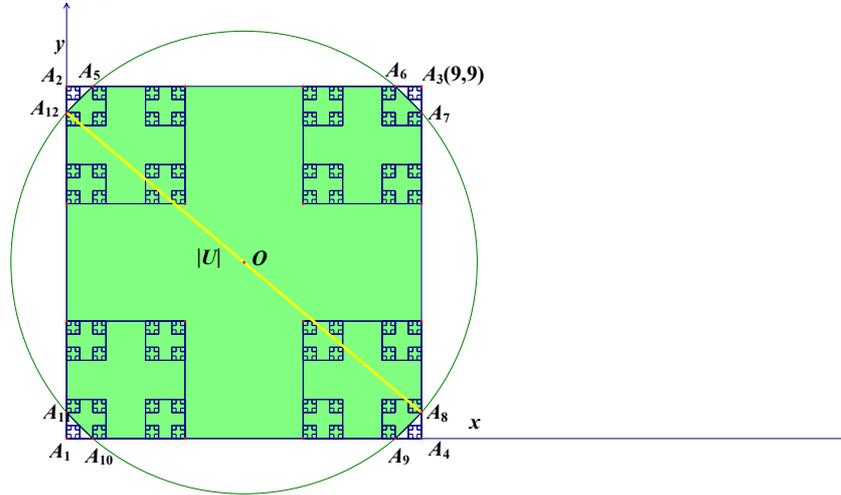

Use this square to construct $C \times C$, pick the middle point of $A_1A_3$, O, and pick $A_5$ (2/3，9). Construct a circle with center $O$ radius $OA_5$ intersecting $C \times C$ at points $A_5$、$A_6$、…、$A_{12}$. Take $\odot O \cap C \times C = U$, where $U$ is the cover set.

Obviously the diameter of $U$ is $A_5A_9$, given by

$$|U| = \sqrt{\left(9 - \frac{4}{3}\right)^2 + 9^2} = \frac{\sqrt{1258}}{3}, r = \frac{|U|}{2} = \frac{\sqrt{1258}}{6}.$$

So the equation of $\odot O$ is

$$\left(x - \frac{9}{2}\right)^2 + \left(y - \frac{9}{2}\right)^2 = \frac{629}{18}.$$

According to the symmetry of $C \times C$, we only need to measure the ratio of squares that *do not contain* in set $U$ on the top left corner, then times 4 and subtracted it from 1.

**The 1st level square in here is the square $A_1A_2A_3A_4$ with side length 9.**

Obviously, some squares on 4th level are not fully contained by $U$.

Taking

$$x = \frac{2}{9},$$

substituting $x$ into the equation of $\odot O$ and picking the larger $y$, we have

$$y = \frac{9}{2} + \frac{\sqrt{5393}}{18} \approx 8.5798 < 9 - \frac{2}{9} \approx 8.7778.$$

This result claims that $U$ dose not fully cover the 5th level small squares on the bottom right corner of the 4th level small square on the top left corner.

Taking

$$x = \frac{2}{9},$$

substituting $x$ into the equation of $\odot O$ and picking the larger $y$, we have

$$y = \frac{9}{2} + \frac{\sqrt{5393}}{18} \approx 8.5798 < 9 - \frac{1}{3} \approx 8.6667.$$

Taking
$$y = 9 - \frac{1}{3} = \frac{26}{3},$$
substituting $y$ into the equation of $\odot O$ and picking smaller $y$, we have
$$x = \frac{9}{2} - \frac{\sqrt{633}}{6} \approx 0.3068 > \frac{2}{9} + \frac{2}{27} \approx 0.2963.$$

By the symmetry of $\odot O$, this result claims that $U$ dose not fully cover the 6$^{th}$ level small squares on the bottom right corner of the 4$^{th}$ level small square on the top left corner, as shown in the figure below.

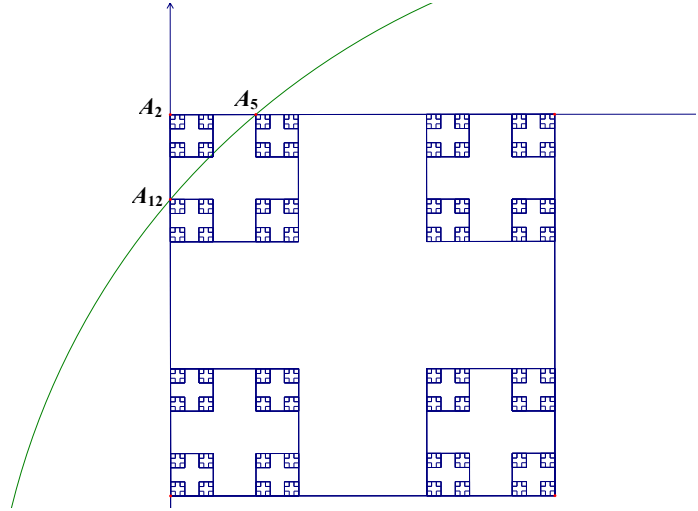

As a result, the cover of 4$^{th}$ level small squares on the top left corner by $U$ should be considered to 7$^{th}$ level squares.

Taking
$$x = \frac{2}{9} + \frac{2}{27} = \frac{8}{27},$$
substituting $x$ into the equation of $\odot O$ and picking the larger $y$, we have
$$y = \frac{9}{2} + \frac{\sqrt{50369}}{54} \approx 8.6561 < 9 - \frac{1}{3} \approx 8.6667.$$

Taking
$$x = \frac{2}{9} + \frac{2}{27} + \frac{2}{81} = \frac{26}{81},$$
substituting $x$ into the equation of $\odot O$ and picking the larger $y$, we have
$$y = \frac{9}{2} + \frac{\sqrt{458753}}{162} \approx 8.6809 > 9 - \frac{2}{9} - \frac{2}{27} - \frac{2}{81} \approx 8.6790.$$

Taking
$$x = \frac{2}{9} + \frac{2}{27} + \frac{1}{81} = \frac{25}{81},$$
substituting $x$ into the equation of $\odot O$ and picking the larger $y$, we have
$$y = \frac{9}{2} + \frac{\sqrt{456041}}{162} \approx 8.6686 < 9 - \frac{2}{9} - \frac{2}{27} - \frac{1}{81} \approx 8.6914.$$

This result claims that, $U$ fully covers the $7^{th}$ level small square on the bottom right corner of the $4^{th}$ level small square on the top left corner, besides $U$ also cover part of two $7^{th}$ level small squares adjacent to this $7^{th}$ level small square, as shown in the figure below.

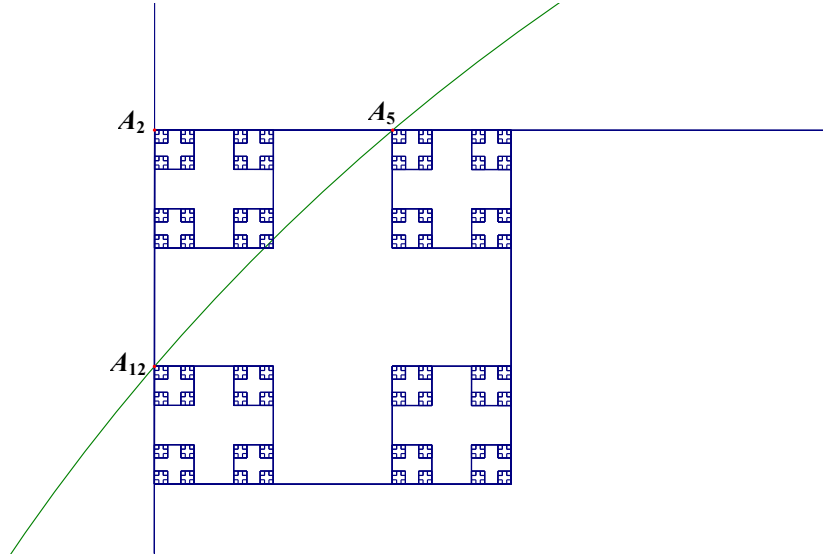

Consider the coverage of $8^{th}$ level squares.

Because of the symmetry of $\odot O$, consider the coverage of the $8^{th}$ level small square in the $7^{th}$ level small square on the left side of the $7^{th}$ level small square on the bottom right corner of the $4^{th}$ level small square on the top left corner.

The calculation above shows that, when
$$y = 9 - \frac{1}{3} = \frac{26}{3},$$
the value of $x$ is
$$x = \frac{9}{2} - \frac{\sqrt{633}}{6} \approx 0.3068 > \frac{2}{9} + \frac{2}{27} + \frac{1}{243} \approx 0.3004,$$
$$x = \frac{9}{2} - \frac{\sqrt{633}}{6} \approx 0.3068 > \frac{2}{9} + \frac{2}{27} + \frac{2}{243} \approx 0.3045.$$
When
$$x = \frac{2}{9} + \frac{2}{27} + \frac{1}{81} = \frac{25}{81},$$
the largest value of $y$ is
$$y = \frac{9}{2} + \frac{\sqrt{456041}}{162} \approx 8.6686 < 9 - \frac{2}{9} - \frac{2}{27} - \frac{2}{81} - \frac{2}{243} \approx 8.6708.$$
Taking
$$x = \frac{2}{9} + \frac{2}{27} + \frac{1}{243} = \frac{73}{243},$$
substituting $x$ into the equation of $\odot O$ and picking the larger $y$, we have
$$y = \frac{9}{2} + \frac{\sqrt{4088057}}{486} \approx 8.6602 < 9 - \frac{2}{9} - \frac{2}{27} - \frac{2}{81} - \frac{1}{243} \approx 8.6749.$$
This result shows that these $8^{th}$ level small squares are all covered part of itself by $U$, as shown in the figure below.

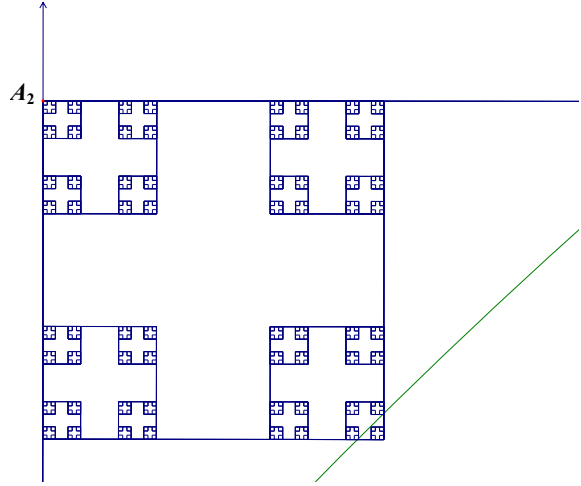

Consider the coverage of 9th level squares.
The calculation above shows that, when
$$x = \frac{2}{9} + \frac{2}{27} + \frac{1}{81} = \frac{25}{81},$$
the value of $y$ is
$$y = \frac{9}{2} + \frac{\sqrt{456041}}{162} \approx 8.6686 < 9 - \frac{2}{9} - \frac{2}{27} - \frac{2}{81} - \frac{2}{243} - \frac{1}{729} \approx 8.6694,$$
$$y = \frac{9}{2} + \frac{\sqrt{456041}}{162} \approx 8.6686 > 9 - \frac{2}{9} - \frac{2}{27} - \frac{2}{81} - \frac{2}{243} - \frac{2}{729} \approx 8.6680.$$
Taking
$$x = \frac{2}{9} + \frac{2}{27} + \frac{2}{243} + \frac{2}{729} = \frac{224}{729},$$
substituting $x$ into the equation of $\odot O$ and picking the larger $y$, we have
$$y = \frac{9}{2} + \frac{\sqrt{36914873}}{1458} \approx 8.6672 < 9 - \frac{2}{9} - \frac{2}{27} - \frac{2}{81} - \frac{2}{243} - \frac{2}{729} \approx 8.6680.$$
According to the symmetry of $\odot O$, this result claims that there are 8 9th level small squares partly covered by set $U$ in the whole $C \times C$, these 9th level small squares are mainly located in the 7th level small square on the bottom right corner of the 4th level small square on the top left corner and two 7th small squares near this 7th level small square Consider the coverage of 10th level squares.
Taking
$$x = \frac{2}{9} + \frac{2}{27} + \frac{2}{243} + \frac{2}{729} + \frac{1}{2187} = \frac{673}{2187},$$
substituting $x$ into the equation of $\odot O$ and picking the larger $y$, we have
$$y = \frac{9}{2} + \frac{\sqrt{332307209}}{4374} \approx 8.6677 > 9 - \frac{2}{9} - \frac{2}{27} - \frac{2}{81} - \frac{2}{243} - \frac{2}{729} - \frac{1}{2187}$$
$$\approx 8.6676.$$
Taking
$$x = \frac{2}{9} + \frac{2}{27} + \frac{2}{243} + \frac{2}{729} = \frac{224}{729},$$
substituting $x$ into the equation of $\odot O$ and picking the larger $y$, we have

$$y = \frac{9}{2} + \frac{\sqrt{36914873}}{1458} \approx 8.6672 < 9 - \frac{2}{9} - \frac{2}{27} - \frac{2}{81} - \frac{2}{243} - \frac{2}{729} - \frac{1}{2187}$$
$$\approx 8.6676,$$
$$y = \frac{9}{2} + \frac{\sqrt{36914873}}{1458} \approx 8.6672 > 9 - \frac{2}{9} - \frac{2}{27} - \frac{2}{81} - \frac{2}{243} - \frac{2}{729} - \frac{2}{2187}$$
$$\approx 8.6671.$$

Taking
$$x = \frac{2}{9} + \frac{2}{27} + \frac{2}{243} + \frac{2}{729} + \frac{2}{2187} = \frac{674}{2187},$$
substituting $x$ into the equation of $\odot O$ and picking the larger $y$, we have
$$y = \frac{9}{2} + \frac{\sqrt{332380553}}{4374} \approx 8.66806 > 9 - \frac{2}{9} - \frac{2}{27} - \frac{2}{81} - \frac{2}{243} - \frac{2}{729} \approx 8.66804.$$

This result clearly points out that, there are three $10^{th}$ level small squares coved by $U$ lie on the top left corner of $C \times C$, and they are located in the $7^{th}$ level small square on the bottom right corner of the $4^{th}$ level small square on the top left corner and two $7^{th}$ small squares near this $7^{th}$ level small square, respectively. In other word, $U$ contains $4^3 + 3 \times 2 = 70$ $10^{th}$ level small squares in the $6^{th}$ level small square in the bottom right corner of the $4^{th}$ level small square on the top left corner; that is, the total number of $10^{th}$ level small squares lying on the top left corner and outside of set $U$ is $(16 - 1) \times 4^4 + 4^4 - 70 = 186 + 3840 = 4026$. In this case, the total number of $10^{th}$ level small squares lying outside set $U$ is $4026 \times 4 = 16104$.

Notice that the total number of $10^{th}$ level small squares is
$$N_{10} = 4^{10-1} = 4^9 = 262144.$$
So, the fraction of $10^{th}$ level small squares lying inside set $U$ is
$$1 - \frac{16104}{262144} = \frac{30755}{32768}.$$

Then, the diameter of "$U$" of the original $C \times C$ is (divided the diameter of $U$ of square $A_1A_2A_3A_4$ by similarity ratio 9)
$$|U| = \frac{\sqrt{\left(9 - \frac{4}{3}\right)^2 + 9^2}}{9} = \frac{\sqrt{1258}}{27}.$$

Then by Lemma 3, we have
$$\frac{30755}{32768} H^s(C \times C) \leq |U|^s = \left(\frac{\sqrt{1258}}{27}\right)^{\log_3 4}.$$

Solving this inequity, we have
$$\boldsymbol{H^s(C \times C) \leq 1.503263}.$$

In this paper, we consider the $10^{th}$ level small squares. If we consider smaller squares, the upper bound will be smaller than this result!

**6.5 Theorem 4**

Combine the ideas of method 1 and method 2, use infinite series to make sure the diameter and construct a circle in the same way of method 2. The Hausdorff measure of the Cartesian product of the middle third Cantor set with itself satisfies the inequality
$$\boldsymbol{H^s(C \times C) \leq 1.502483}.$$

**Proof:** Firstly, construct the rectangular coordinator system $xA_1y$ with point $A_1(0,0)$, then construct square $A_1A_2A_3A_4$ with $A_1A_2=1$ by point $A_1$ for vertex and $x$, $y$ axes for sides, as shown in the figure below.

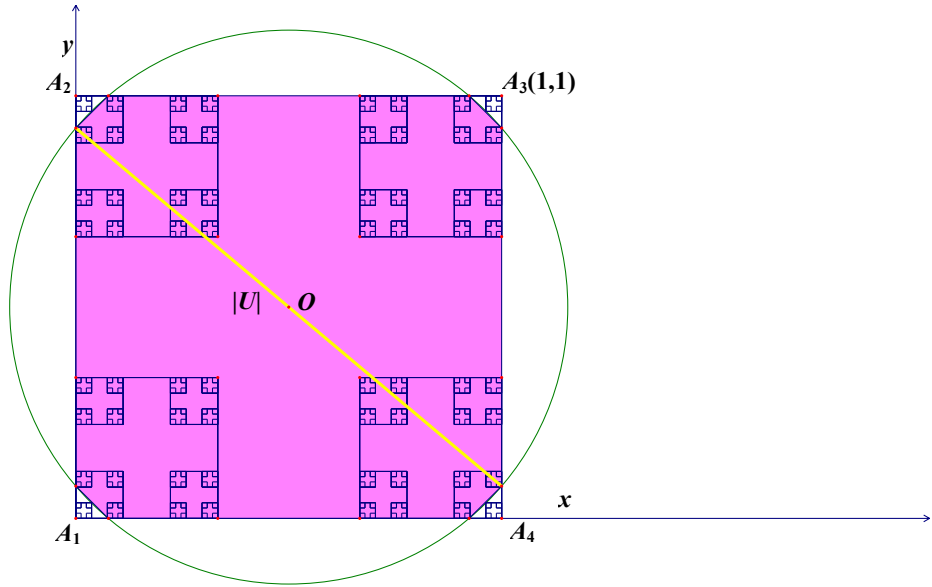

According to method 1 and method 2, we take

$$\varphi = \sum_{i=1}^{\infty} \frac{1}{3^{3i}} = \frac{1}{26},$$

chose the middle point of $A_1A_3$ as $O$, and construct a circle with center $O$ and radius

$$r = \frac{\sqrt{(1-4\varphi)^2 + 1^2}}{2} = \frac{\sqrt{290}}{26}.$$

So, the equation of this circle is

$$\left(x - \frac{1}{2}\right)^2 + \left(y - \frac{1}{2}\right)^2 = \frac{290}{676}.$$

According to the symmetry of $U$, we only need to consider the coverage of small squares on the top left corner. Because the diameter of $U$ is determined by 4[th] level small squares, we first consider the coverage of 5[th] small squares. **The 1[st] level square here is the square $A_1A_2A_3A_4$ with side length 1.**
Taking

$$x = \frac{2}{27},$$

substituting $x$ into the equation of $\odot O$ and picking the larger $y$, we have

$$y = \frac{1}{2} + \frac{\sqrt{122009}}{702} \approx 0.9976 > 1 - \frac{1}{81} \approx 0.9877.$$

Taking

$$x = \frac{1}{27},$$

substituting $x$ into the equation of $\odot O$ and picking the larger $y$, we have

$$y = \frac{1}{2} + \frac{\sqrt{105785}}{702} \approx 0.9633 > 1 - \frac{1}{27} \approx 0.9630,$$

$$y = \frac{1}{2} + \frac{\sqrt{105785}}{702} \approx 0.9633 < 1 - \frac{1}{27} + \frac{1}{81} \approx 0.9753.$$

When

$$y = 1,$$

compare $x$ with $2\varphi$

$$2\varphi = \frac{2}{26} \approx 0.07692.$$

Take

$$x = \frac{2}{27} + \frac{1}{81} \approx 0.08642 > 2\varphi.$$

According to the symmetry of this figure, we know that $\odot O$ must intersect with three $4^{th}$ level small squares lying on the top left corner, bottom left corner and top right corner of the $3^{rd}$ level small square lying on the top left corner, respectively. The intersections lie on the side of the $5^{th}$ level small square on the bottom right corner of the $4^{th}$ level small square on the top left corner and the sides of the small squares in another two $4^{th}$ level small squares lie beside the $4^{th}$ level small squares on the top left corner, as shown in the figure below.

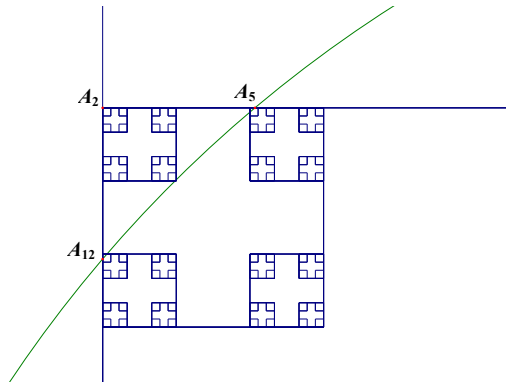

Consider the coverage of $6^{th}$ level squares.
The calculation above shows that, when

$$x = \frac{1}{27},$$

the value of $y$ is

$$y = \frac{1}{2} + \frac{\sqrt{105785}}{702} \approx 0.9633 < 1 - \frac{1}{27} + \frac{1}{243} \approx 0.9671.$$

When

$$x = \frac{2}{27},$$

the value of $y$ is

$$y = \frac{1}{2} + \frac{\sqrt{122009}}{702} \approx 0.9976 > 1 - \frac{1}{243} \approx 0.9959.$$

When

$$y = 1,$$

compare $x$ with $2\varphi$

$$2\varphi = \frac{2}{26} \approx 0.07692.$$

Take

$$x = \frac{2}{27} + \frac{1}{243} \approx 0.07819 > 2\varphi.$$

According to the symmetry of the figure, we know that $\odot O$ passes through three 6$^{th}$ level small squares in these 5$^{th}$ level small squares claimed above, that is, this circle passes through the 6$^{th}$ level small square on the bottom right corner of the 4$^{th}$ level small square on the top left corner and two 6$^{th}$ level small squares on the top left corner of two 4$^{th}$ level small squares lie beside the 4$^{th}$ level small square on the top left corner, respectively, as shown in the figure below.

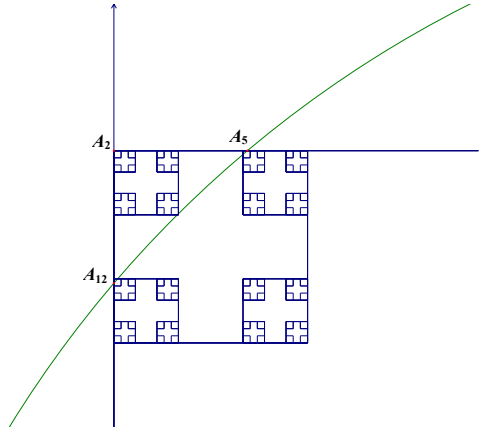

Consider the coverage of 7$^{th}$ level squares.
The calculation above shows that, when

$$x = \frac{1}{27},$$

the value of $y$ is

$$y = \frac{1}{2} + \frac{\sqrt{105785}}{702} \approx 0.9633 < 1 - \frac{1}{27} + \frac{1}{729} \approx 0.9643.$$

When

$$x = \frac{2}{27},$$

the value of $y$ is

$$y = \frac{1}{2} + \frac{\sqrt{122009}}{702} \approx 0.9976 < 1 - \frac{1}{729} \approx 0.9986,$$

$$y = \frac{1}{2} + \frac{\sqrt{122009}}{702} \approx 0.9976 > 1 - \frac{2}{729} \approx 0.9973.$$

When

$$y = 1,$$

compare $x$ with $2\varphi$

$$2\varphi = \frac{2}{26} \approx 0.07692.$$

Taking

$$x = \frac{2}{27} + \frac{1}{729} \approx 0.07545 < 2\varphi,$$

$$x = \frac{2}{27} + \frac{2}{729} \approx 0.07682 < 2\varphi.$$

Taking

$$x = \frac{2}{27} + \frac{1}{729} = \frac{55}{729},$$

substituting $x$ into the equation of $\odot O$ and picking the larger $y$, we have

$$y = \frac{1}{2} + \frac{\sqrt{89363681}}{18954} \approx 0.9987 > 1 - \frac{1}{729} \approx 0.9986.$$

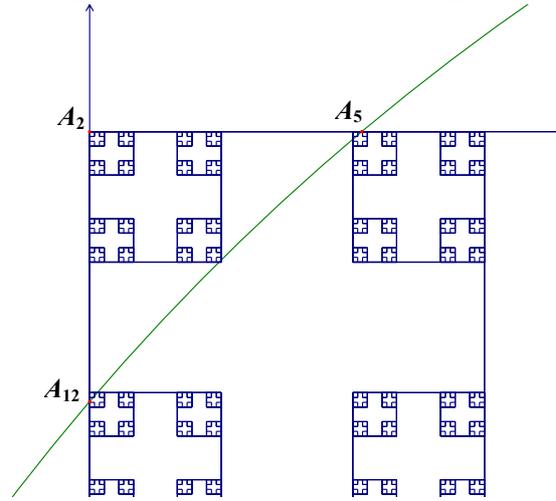

According to the symmetry of the figure, we know that $\odot O$ passes through three 7th level small squares in these 5th level small squares claimed above and part of 7th level small squares, respectively. Also, $\odot O$ passes through the 7th level small square on the bottom right corner of the 4th level small square on the top left corner, as shown in the figure above.

Consider the coverage of 8th level squares.

Taking

$$x = \frac{2}{27} + \frac{1}{2187} = \frac{55}{729},$$

substituting $x$ into the equation of $\odot O$ and picking the larger $y$, we have

$$y = \frac{1}{2} + \frac{\sqrt{801759761}}{56862} \approx 0.9980 < 1 - \frac{1}{729} \approx 0.9986.$$

Taking

$$x = \frac{2}{27} + \frac{2}{729} = \frac{56}{729},$$

substituting $x$ into the equation of $\odot O$ and picking the larger $y$, we have

$$y = \frac{1}{2} + \frac{\sqrt{89781449}}{18954} \approx 0.9999 > 1 - \frac{1}{729} \approx 0.9995.$$

Taking

$$y = 1 - \frac{1}{729} = \frac{728}{729},$$

substituting $y$ into the equation of $\odot O$ and picking the smaller $x$, we have

$$x = \frac{1}{2} - \frac{\sqrt{64796489}}{18954} \approx 0.07531 < \frac{2}{27} + \frac{1}{729} \approx 0.07545,$$

$$x = \frac{1}{2} - \frac{\sqrt{64796489}}{18954} \approx 0.07531 > \frac{2}{27} + \frac{1}{2187} \approx 0.07453,$$

$$x = \frac{1}{2} - \frac{\sqrt{64796489}}{18954} \approx 0.07531 > \frac{2}{27} + \frac{2}{2187} \approx 0.07499.$$

The calculation above shows that, when

$$x = \frac{1}{27},$$

the value of $y$ is

$$y = \frac{1}{2} + \frac{\sqrt{105785}}{702} \approx 0.9633 < 1 - \frac{1}{27} + \frac{1}{2187} \approx 0.9634.$$

When

$$y = 1,$$

compare $x$ with $2\varphi$

$$2\varphi = \frac{2}{26} \approx 0.07692.$$

Take

$$x = \frac{2}{27} + \frac{2}{729} \approx 0.07682 < 2\varphi,$$

$$x = \frac{2}{27} + \frac{2}{729} + \frac{1}{2187} \approx 0.07728 > 2\varphi.$$

According to the symmetry, we know that $\odot O$ passes through 8th level small squares on the bottom right corner of these two 7th level small squares claimed above and the 8th level small squares lie beside these two 7th level small squares, respectively. Also, $\odot O$ passes through the 8th level small square on the bottom right corner of the 4th level small square on the top left corner, as shown in the figure below.

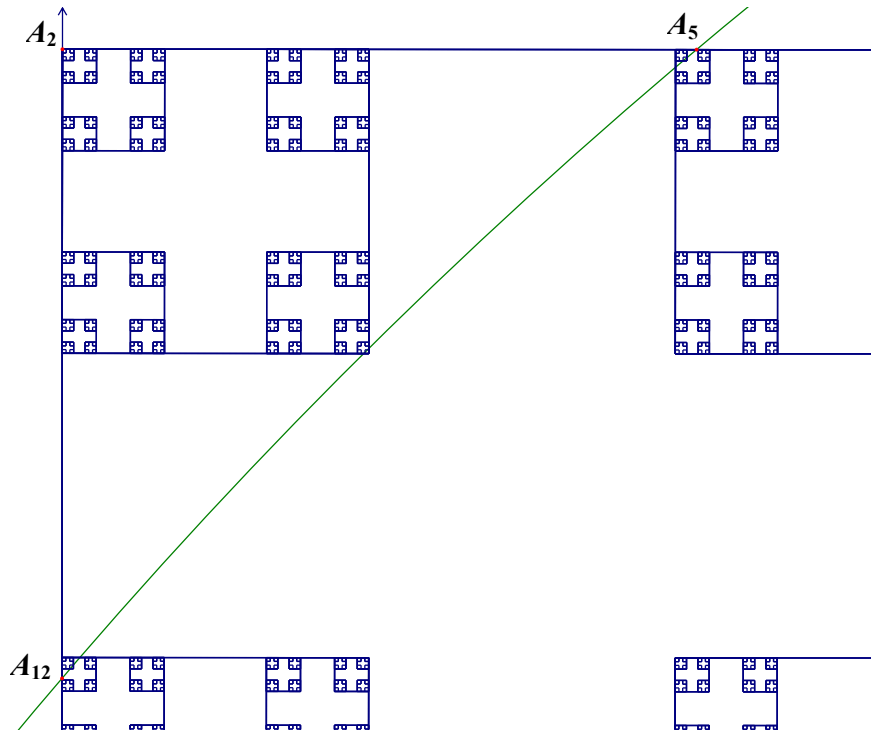

Consider the coverage of $9^{th}$ level squares.
Taking
$$x = \frac{2}{27} + \frac{2}{729} + \frac{1}{6561} = \frac{505}{6561},$$
substituting $x$ into the equation of $\odot O$ and picking the larger $y$, we have
$$y = \frac{1}{2} + \frac{\sqrt{7276050521}}{170586} \approx 1.00004 > 1.$$
Taking
$$x = \frac{2}{27} + \frac{2}{2187} + \frac{2}{6561} = \frac{494}{6561},$$
substituting $x$ into the equation of $\odot O$ and picking the larger $y$, we have
$$y = \frac{1}{2} + \frac{\sqrt{7234691489}}{170586} \approx 0.9986 < 1 - \frac{2}{2187} - \frac{2}{6561} \approx 0.9988.$$
The calculation above shows that, when
$$y = 1 - \frac{1}{729} = \frac{728}{729},$$
the value of $x$ is
$$x = \frac{1}{2} - \frac{\sqrt{64796489}}{18954} \approx 0.07531 > \frac{2}{27} + \frac{2}{2187} + \frac{2}{6561} \approx 0.07529,$$
$$x = \frac{1}{2} - \frac{\sqrt{64796489}}{18954} \approx 0.07531 < \frac{2}{27} + \frac{1}{729} \approx 0.07545.$$
The calculation above shows that, when
$$x = \frac{2}{27} + \frac{1}{729} = \frac{55}{729},$$
the value of $y$ is
$$y = \frac{1}{2} + \frac{\sqrt{89363681}}{18954} \approx 0.9987 > 1 - \frac{1}{729} \approx 0.9986,$$
$$y = \frac{1}{2} + \frac{\sqrt{89363681}}{18954} \approx 0.9987 < 1 - \frac{2}{2187} - \frac{2}{6561} \approx 0.9988.$$
The calculation above shows that, when
$$x = \frac{2}{27} + \frac{2}{729} = \frac{56}{729},$$
the value of $y$ is
$$y = \frac{1}{2} + \frac{\sqrt{89781449}}{18954} \approx 0.9999 > 1 - \frac{1}{6561} \approx 0.9998.$$
When
$$y = 1,$$
compare $x$ with $2\varphi$
$$2\varphi = \frac{2}{26} \approx 0.07692.$$
Take
$$x = \frac{2}{27} + \frac{2}{729} + \frac{1}{6561} \approx 0.07697 < 2\varphi.$$
The calculation above shows that, when

$$x = \frac{1}{27},$$

the value of $y$ is

$$y = \frac{1}{2} + \frac{\sqrt{105785}}{702} \approx 0.9633 < 1 - \frac{1}{27} + \frac{1}{2187} \approx 0.9634.$$

Now, it is not meaningful to subdivide the squares, because the arc is approximately coincided with the diagonals of two 9$^{th}$ level squares and one 8$^{th}$ level small square. Also, this arc intersects with another two 9$^{th}$ level small squares, but the intersection area is extremely small.

For arbitrary positive integer $n$, let the number of $n^{th}$ level small squares be $N$, then the relationship between $N$ and $n$ is

$$N_n = 4^{n-1}.$$

As a result, the total number of 9$^{th}$ level small squares is

$$N_{10} = 4^{9-1} = 4^8 = 65536.$$

According to the symmetry of $U$, we first consider the un-coverage of $C \times C$ on the top right corner. Notice that the arc approximately coincides with the diagonals of two 9$^{th}$ level squares and one 8$^{th}$ level small square, so we approximately pick half of each squares, and get three 9$^{th}$ level small squares. We don't calculate the squares that are half-divided for now. According to the calculation above, there are 15 9$^{th}$ level small squares in the 4$^{th}$ level small square lying on the left of the 4$^{th}$ level small square lying on the top left corner not covered by $U$. According to the symmetry of this figure, we know that there are 15 9$^{th}$ level small squares in the 4$^{th}$ level small square lying on under the 4$^{th}$ level small square lying on the top left corner not covered by $U$ too. For the 4$^{th}$ level small square on the top left corner, remove the half 8$^{th}$ level small square. There are

$$4^4 - 1 = 255$$

8$^{th}$ level small squares. So, there are $255 \times 4 = 1020$ 9$^{th}$ level small squares. So, the total number of 9$^{th}$ level small squares that are not covered by $U$ on the top left corner is $1020 + 15 \times 2 + 3 = 1053$. As a result, the total number of 9$^{th}$ level small squares that are not covered by $U$ in the whole $C \times C$ is $1053 \times 4 = 4212$. So, the fraction of 9$^{th}$ level small squares covered by $U$ is

$$1 - \frac{4212}{65536} = \frac{15331}{16384},$$

and the diameter of $U$ is

$$|U| = 2r = \frac{\sqrt{290}}{13}.$$

Then according to Lemma 3, we have

$$\frac{15331}{16384} H^s(C \times C) \leq |U|^s = \left(\frac{\sqrt{290}}{13}\right)^{\log_3 4}.$$

Solving this inequity, we have

$$\boldsymbol{H^s(C \times C) \leq 1.502483}.$$

The upper bound of

$$H^s(C \times C)$$

in this paper is slightly different from the precise upper bound of this method, because these 9[th] level small squares are not exactly half-divided, and two 9[th] level small squares are not fully covered (In this article, we take the fully cover, because the uncovered area is so small that it won't influence the result too much).

**6.6 Main Theorem**

In conclusion, we use four different methods (constructs octagon, takes infinity, constructs circle, and expands the coverage but doesn't change the diameter, sums up method 2 and method 3 to create method 4). As a result, under different method, the author gets 4 upper bounds:

**The Hausdorff measure of the Cartesian product of the middle third Cantor set with itself satisfies the inequalities**

$$\begin{cases} H^s(C \times C) \leq 1.548563 \\ H^s(C \times C) \leq 1.504975 \\ H^s(C \times C) \leq 1.502878 \\ H^s(C \times C) \leq 1.503263 \\ H^s(C \times C) \leq 1.502483 \end{cases}$$

Take the most optimized result, the main theorem of this paper is the following

The Hausdorff measure of the Cartesian product of the middle third Cantor set with itself satisfies the inequality

$$H^s(C \times C) \leq 1.502483.$$

# 7 Conclusion

we use three different methods, and gets five different upper bounds

**The Hausdorff measure of the Cartesian product of the middle third Cantor set with itself satisfies the inequalities**

$$\begin{cases} H^s(C \times C) \leq 1.548563 \\ H^s(C \times C) \leq 1.504975 \\ H^s(C \times C) \leq 1.502878 \\ H^s(C \times C) \leq 1.503263 \\ H^s(C \times C) \leq 1.502483 \end{cases}$$

The description of these five different upper bounds:

1. After the correction of an inaccurate upper bound, use a "basic interval" of $n^{th}$ level $C \times C$ as the cover set, where $n \in N^+$, we have:

The Hausdorff measure of the Cartesian product of the middle third Cantor set with itself satisfies the inequality

$$H^s(C \times C) \leq 1.548563.$$

2. Use octagon as the cover set *U*, we have

The Hausdorff measure of the Cartesian product of the middle third Cantor set with itself satisfies the inequality

$$H^s(C \times C) \leq 1.504975.$$

3. Use octagon and isometric series which is used to ensure the diameter (the methods are shown in the paper), we obtaining the following:

The Hausdorff measure of the Cartesian product of the middle third Cantor set with itself satisfies the inequality

$$H^s(C \times C) \leq 1.502878.$$

4. Keep the diameter constantly, and take circle as the bigger cover set, because circle has the biggest area compared to other sets having same diameter. Using this method, we obtain the following

The Hausdorff measure of the Cartesian product of the middle third Cantor set with itself satisfies the inequality

$$H^s(C \times C) \leq 1.503263.$$

5. Combine method 1 and 2 together and take both circle and isometric series.

The Hausdorff measure of the Cartesian product of the middle third Cantor set with itself satisfies the inequality

$$H^s(C \times C) \leq 1.502483.$$

Taking the best upper bound, we have

**The Hausdorff measure of the Cartesian product of the middle third Cantor set with itself satisfies the inequality**

$$H^s(C \times C) \leq 1.502483.$$

As far as we know, this upper bound is the smallest upper bound.

The creation of this paper is, **how to take the cover set *U*, so that the upper bound can be smaller**. In addition, this paper also corrects an inaccurate upper bound of a reference. These points show that this paper makes some contributions to the advanced mathematical topics.

## 8 Postscript

**Author**
Fan Yuchen, a 16-year-old grade 11 student from International Department of the Affiliated High School of South China Normal University, graduated from The Middle School Attached to Sun Yat-sen University. Mailbox fanyc.frank2017@outlook.com.
**Instructor**
Baoguo Jia, professor of mathematics in Sun Yat-sen University. Major research field is Geometry and topology. Mailbox is mcsjbg@mail.sysu.edu.cn

# 10 Appendix

**Matlab code**

```
>> syms k;y=(sqrt(2*3^(2*k)-12*3^k+26)/(3^k-1))^(log(4)/log(3))*(4^k-2)/(4^k-6)

y =

((2*3^(2*k)-12*3^k+26)^(1/2)/(3^k-1))^(2841455003081375/2251799813685248)*(4^k-2)/(4^k-6)

>> diff(y)

ans =

2841455003081375/2251799813685248*((2*3^(2*k)-12*3^k+26)^(1/2)/(3^k-1))^(589655189396127/2251799813685248)*(4^k-2)/(4^k-6)*(1/2/(2*3^(2*k)-12*3^k+26)^(1/2)/(3^k-1)*(4*3^(2*k)*log(3)-12*3^k*log(3))-(2*3^(2*k)-12*3^k+26)^(1/2)/(3^k-1)^2*3^k*log(3))+((2*3^(2*k)-12*3^k+26)^(1/2)/(3^k-1))^(2841455003081375/2251799813685248)*4^k*log(4)/(4^k-6)-((2*3^(2*k)-12*3^k+26)^(1/2)/(3^k-1))^(2841455003081375/2251799813685248)*(4^k-2)/(4^k-6)^2*4^k*log(4)

>>                k=solve('2841455003081375/2251799813685248*((2*3^(2*k)-12*3^k+26)^(1/2)/(3^k-1))^(589655189396127/2251799813685248)*(4^k-2)/(4^k-6)*(1/2/(2*3^(2*k)-12*3^k+26)^(1/2)/(3^k-1)*(4*3^(2*k)*log(3)-12*3^k*log(3))-(2*3^(2*k)-12*3^k+26)^(1/2)/(3^k-1)^2*3^k*log(3))+((2*3^(2*k)-12*3^k+26)^(1/2)/(3^k-1))^(2841455003081375/2251799813685248)*4^k*log(4)/(4^k-6)-((2*3^(2*k)-12*3^k+26)^(1/2)/(3^k-1))^(2841455003081375/2251799813685248)*(4^k-2)/(4^k-6)^2*4^k*log(4)=0','k')

k =

2.7805145062820030104236744396445

>> diff(diff(y))

ans =
```

16754786880025208040104908346325/5070602400912917605986812821504/((2*3^(2*k)-12*3^k+26)^(1/2)/(3^k-1))^(1662144624289121/2251799813685248)*(4^k-2)/(4^k-6)*(1/2/(2*3^(2*k)-12*3^k+26)^(1/2)/(3^k-1)*(4*3^(2*k)*log(3)-12*3^k*log(3))-(2*3^(2*k)-12*3^k+26)^(1/2)/(3^k-1)^2*3^k*log(3))^2+2841455003081375/1125899906842624*((2*3^(2*k)-12*3^k+26)^(1/2)/(3^k-1))^(589655189396127/2251799813685248)*4^k*log(4)/(4^k-6)*(1/2/(2*3^(2*k)-12*3^k+26)^(1/2)/(3^k-1)*(4*3^(2*k)*log(3)-12*3^k*log(3))-(2*3^(2*k)-12*3^k+26)^(1/2)/(3^k-1)^2*3^k*log(3))-2841455003081375/1125899906842624*((2*3^(2*k)-12*3^k+26)^(1/2)/(3^k-1))^(589655189396127/2251799813685248)*(4^k-2)/(4^k-6)^2*(1/2/(2*3^(2*k)-12*3^k+26)^(1/2)/(3^k-1)*(4*3^(2*k)*log(3)-12*3^k*log(3))-(2*3^(2*k)-12*3^k+26)^(1/2)/(3^k-1)^2*3^k*log(3))*4^k*log(4)+2841455003081375/2251799813685248*((2*3^(2*k)-12*3^k+26)^(1/2)/(3^k-1))^(589655189396127/2251799813685248)*(4^k-2)/(4^k-6)*(-1/4/(2*3^(2*k)-12*3^k+26)^(3/2)/(3^k-1)*(4*3^(2*k)*log(3)-12*3^k*log(3))^2-1/(2*3^(2*k)-12*3^k+26)^(1/2)/(3^k-1)^2*(4*3^(2*k)*log(3)-12*3^k*log(3))*3^k*log(3)+1/2/(2*3^(2*k)-12*3^k+26)^(1/2)/(3^k-1)*(8*3^(2*k)*log(3)^2-12*3^k*log(3)^2)+2*(2*3^(2*k)-12*3^k+26)^(1/2)/(3^k-1)^3*(3^k)^2*log(3)^2-(2*3^(2*k)-12*3^k+26)^(1/2)/(3^k-1)^2*3^k*log(3)^2)+((2*3^(2*k)-12*3^k+26)^(1/2)/(3^k-1))^(2841455003081375/2251799813685248)*4^k*log(4)^2/(4^k-6)-2*((2*3^(2*k)-12*3^k+26)^(1/2)/(3^k-1))^(2841455003081375/2251799813685248)*(4^k)^2*log(4)^2/(4^k-6)^2+2*((2*3^(2*k)-12*3^k+26)^(1/2)/(3^k-1))^(2841455003081375/2251799813685248)*(4^k-2)/(4^k-6)^3*(4^k)^2*log(4)^2-((2*3^(2*k)-12*3^k+26)^(1/2)/(3^k-1))^(2841455003081375/2251799813685248)*(4^k-2)/(4^k-6)^2*4^k*log(4)^2

>> z=ans

z =

16754786880025208040104908346325/5070602400912917605986812821504/((2*3^(2*k)-12*3^k+26)^(1/2)/(3^k-1))^(1662144624289121/2251799813685248)*(4^k-2)/(4^k-6)*(1/2/(2*3^(2*k)-12*3^k+26)^(1/2)/(3^k-1)*(4*3^(2*k)*log(3)-12*3^k*log(3))-(2*3^(2*k)-12*3^k+26)^(1/2)/(3^k-1)^2*3^k*log(3))^2+2841455003081375/1125899906842624*((2*3^(2*k)-12*3^k+26)^(1/2)/(3^k-1))^(589655189396127/2251799813685248)*4^k*log(4)/(4^k-6)*(1/2/(2*3^(2*k)-12*3^k+26)^(1/2)/(3^k-1)*(4*3^(2*k)*log(3)-12*3^k*log(3))-(2*3^(2*k)-12*3^k+26)^(1/2)/(3^k-1)^2*3^k*log(3))-2841455003081375/1125899906842624*((2*3^(2*k)-12*3^k+26)^(1/2)/(3^k-1))^(589655189396127/2251799813685248)*(4^k-2)/(4^k-6)^2*(1/2/(2*3^(2*k)-12*3^k+26)^(1/2)/(3^k-1)*(4*3^(2*k)*log(3)-12*3^k*log(3))-(2*3^(2*k)-12*3^k+26)^(1/2)/(3^k-

1)^2*3^k*log(3))*4^k*log(4)+2841455003081375/2251799813685248*((2*3^(2*k)-12*3^k+26)^(1/2)/(3^k-1))^(589655189396127/2251799813685248)*(4^k-2)/(4^k-6)*(-1/4/(2*3^(2*k)-12*3^k+26)^(3/2)/(3^k-1)*(4*3^(2*k)*log(3)-12*3^k*log(3))^2-1/(2*3^(2*k)-12*3^k+26)^(1/2)/(3^k-1)^2*(4*3^(2*k)*log(3)-12*3^k*log(3))*3^k*log(3)+1/2/(2*3^(2*k)-12*3^k+26)^(1/2)/(3^k-1)*(8*3^(2*k)*log(3)^2-12*3^k*log(3)^2)+2*(2*3^(2*k)-12*3^k+26)^(1/2)/(3^k-1)^3*(3^k)^2*log(3)^2-(2*3^(2*k)-12*3^k+26)^(1/2)/(3^k-1)^2*3^k*log(3)^2)+((2*3^(2*k)-12*3^k+26)^(1/2)/(3^k-1))^(2841455003081375/2251799813685248)*4^k*log(4)^2/(4^k-6)-2*((2*3^(2*k)-12*3^k+26)^(1/2)/(3^k-1))^(2841455003081375/2251799813685248)*(4^k)^2*log(4)^2/(4^k-6)^2+2*((2*3^(2*k)-12*3^k+26)^(1/2)/(3^k-1))^(2841455003081375/2251799813685248)*(4^k-2)/(4^k-6)^3*(4^k)^2*log(4)^2-((2*3^(2*k)-12*3^k+26)^(1/2)/(3^k-1))^(2841455003081375/2251799813685248)*(4^k-2)/(4^k-6)^2*4^k*log(4)^2

>> subs(z,2.78051450628200301042367443964 45)

ans =

    0.1063

>> ezplot(y)

**Licenses**